\documentclass[twoside]{article}
\usepackage{amsfonts,amsmath,amssymb, amsthm, amscd, array,  mathrsfs,mathrsfs}

\usepackage[pdftex]{graphicx}
\usepackage{graphics}
\usepackage{cases}
\usepackage{enumerate}
\usepackage{pstricks}
\usepackage{fancyhdr}
\usepackage{geometry}
\usepackage{setspace}
\newtheorem{definition}{Definition}
\newtheorem{theoreme}{Theorem}

\newtheorem{lemme}{Lemma}

\def\R{\mathbb{R}}

\def\F{\mathbb{F}}
\begin{document}
\title{Exact fuzzy solution of the fuzzy heat-like equations}
\author{Lalla Saadia CHADLI, Atimad HARIR and Said MELLIANI}
\date{{\small Laboratoire de Mod\'elisation et Calcul (LMC), Facult\'e des Sciences et Techniques} \\
{\small Universit\'e Sultan Moulay Slimane, BP 523, B\'eni Mellal, Morocco} }

\fancyhf{}
\addtolength{\headwidth}{\marginparsep}
\addtolength{\headwidth}{\marginparwidth}
\fancyhead{} 
\fancyhead[RO,LE]{\thepage}
\fancyhead[CE]{L. S. Chadli, A. Harir and S. Melliani}
\fancyhead[CO]{Exact fuzzy solution of the fuzzy heat-like equations}
\renewcommand{\headrulewidth}{0.0pt}

\maketitle

\begin{abstract}
	In this paper, the Buckley-Feuring method (BFM) and the variational iteration method (VIM)  are used for find exact fuzzy solution of the fuzzy heat-like equations in one and two dimensions. Several examples are given to show the new theorem of Buckley-Feuring solution and the Seikkala solution. The results obtained in all cases show the reliability and the efficiency of this methods.
\end{abstract}

\section{Introduction}
There are strong and efficient techniques to find approximate solutions for the linear and nonlinear equations, that most of these equations don't have exact solution such as heat-like equations. In mathematics ,in order to solve the model of heat-like equations,we will introduce some imprecise parameters. In this work, our idea is solving heat-like equations with fuzzy parameters via the same strategy as Buckley and Feuring [ref,ref] using Variational Iteration Method (VIM).\\
The VIM proposed by He in \cite{raaaa,raaaaa,raaa,raa} ,is a method of solving linear or nonlinear problems \cite{tddd,td,in} and gives rapidly convergent successive approximations of the exact solution if that last exists.\\
The paper is organized as follows: in Section 2, we call some fundamental results on fuzzy numbers. In Sections 3 and 4, fuzzy
heat-like equations and the VIM are illustrated, respectively. In Section 5, the same strategy as in Buckley-Feuring is presented for
two-dimensional fuzzy heat-like equation. Examples are shown in Section 6 and finally conclusion is given in Section 7.

\section{ Preliminary}
We place a bar over a capital letter to denote a fuzzy number of $\R^{n}$. So, $\overline{U}$, $ \overline{K}$, $\overline{\gamma}$, $\overline{\beta}$ etc. all represent fuzzy numbers of  $ \R^{n} $ for some $n$.\\
We write $\mu_{\overline{U}}(t)$, a number in $[0,1]$, for the membership function of $\overline{ U}$ evaluated at $t \in \R^{n}$.
An $\alpha$-cut of $\overline{U}$ is always a closed and bounded interval that written $\overline{U}[\alpha]$, is defined as
$\displaystyle \{ t : \mu_{\overline{U}}(t) \geq \alpha \}$ for $\alpha \in [0,1]$. We separately specify $\overline{U}[0]$ as the closure of the union of all the  $\overline{U}[\alpha]$ for $0< \alpha \leq 1$

\begin{definition} \cite{tttb}
Let $\F(\R)=\{ \overline{U} : \overline{U}:\R \rightarrow [0,1], \text{ satisfies (i)-(iv) } \}:$
  \begin{itemize}
    \item[\text{(i)}] $\forall  \overline{U} \in \F(\R),\quad \overline{U}$ is normal
    \item[(ii)] $\forall  \overline{U} \in \F(\R),\quad \overline{U}$ is a convex fuzzy set
    \item[(iii)] $\forall  \overline{U} \in \F(\R),\quad \overline{U}$ is upper semi-continuous on  $\R$ and
    \item[(iv)] $\overline{U}[\alpha]$ is a compact set $ \forall \alpha \in [0,1]$.
  \end{itemize}
then  $\F(\R)$ is called fuzzy number space and  $\forall \overline{U} \in \F(\R)$, $\overline{U}$ is called a fuzzy number.
\end{definition}

\begin{definition} \cite{tttb,tttc} 
We represent an arbitrary fuzzy number by an ordered pair of functions  $\overline{U}[\alpha]=[u_{1}(\alpha),u_{2}(\alpha)] , \quad \alpha \in [0,1]$ which satisfy the following requirements:
  \begin{itemize}
    \item[(a)] $u_{1}(\alpha)$ is a bounded left continuous nondecreasing function over $[0,1]$,
    \item[(b)] $u_{2}(\alpha)$ is a bounded left continuous nonincreasing function on $[0,1]$,
    \item[(c)] $u_{1}(\alpha) \leq u_{2}(\alpha)$, \quad $0 \leq \alpha \leq 1$
  \end{itemize}
\end{definition}

\begin{definition}
Let $\overline{U}=(u_{1},u_{2},u_{3}), \ ( u_{1}< u_{2} < u_{3})$. $\overline{U}$ is called triangular fuzzy number with peak (or center) $u_{2}$, left width $u_{2}-u_{1} > 0$ and right width  $u_{3}-u_{2} > 0$, if its membership function has the following form :
\begin{equation*}
  \mu_{\overline{U}}(t) = \begin{cases}
      1-\dfrac{u_{2}-t}{u_{2}-u_{1}}, & u_{1} \leq t \leq u_{2} \\
      1-\dfrac{t-u_{2}}{u_{3}-u_{2}}, & u_{2} \leq t \leq u_{3}
\end{cases}
\end{equation*}
The support of $\overline{U}$ is $[u_{1},u_{3}]$. We will write :
\begin{enumerate}
  \item $ \overline{U} > 0 \  if \  u_{1} > 0 $,
  \item $ \overline{U} \geq 0 \  if \  u_{1} \geq 0 $,
  \item $ \overline{U} < 0 \ if \  u_{3} < 0 $,
\item $ \overline{U} \leq 0 \ if \  u_{3} \leq 0 $.
\end{enumerate}
\end{definition}

\begin{definition}
For arbitrary fuzzy numbers $\overline{U}[\alpha]=[u_{1}(\alpha),u_{2}(\alpha)]$ and $\overline{V}[\alpha]=[v_{1}(\alpha),v_{2}(\alpha)]$ we have algebraic operations bellow :
\begin{enumerate}
  \item $ (\overline{U}+\overline{V})[\alpha] = [u_{1}(\alpha)+ v_{1}(\alpha),u_{2}(\alpha) + v_{2}(\alpha)]$
  \item $ (\overline{U}-\overline{V})[\alpha] = [u_{1}(\alpha)- v_{2}(\alpha) ,u_{2}(\alpha) -v_{1}(\alpha) ]$
  \item 
  \begin{equation*}
       k \overline{U}[\alpha] = \begin{cases}
         [k\, u_1(\alpha), k\, u_2(\alpha)] & \text{$k\geq 0$} \\
         [k\, u_2(\alpha), k\, u_1(\alpha)] & \text{$k<0$}
       \end{cases}
   \end{equation*}
  \item $(\overline{U}.\overline{V})[\alpha]$ = $\{ \min z,\max z \}$ with
   $$ z=\{u_{1}(\alpha).v_{1}(\alpha),u_{1}(\alpha).v_{2}(\alpha),u_{2}(\alpha) .v_{1}(\alpha),u_{2}(\alpha). v_{2}(\alpha)\} $$
   \item if $0\notin[v_{1}(\alpha),v_{2}(\alpha)] $  
   $$\frac{\overline{U}}{\overline{V}}[\alpha] =[(\frac{u_{1}}{v_{1}})(\alpha),(\frac{u_{2}}{v_{2}})(\alpha)] $$
where 
$$
(\frac{u_{1}}{v_{1}})(\alpha)=\min\left\{\frac{u_{1}(\alpha)}{v_{1}(\alpha)},
\frac{u_{1}(\alpha)}{v_{2}(\alpha)},\frac{u_{2}(\alpha)}{v_{1}(\alpha)},\frac{u_{2}(\alpha)}{v_{2}(\alpha)}\right\}
$$
$$
(\frac{u_{2}}{v_{2}})(\alpha)=\max\left\{\frac{u_{1}(\alpha)}{v_{1}(\alpha)},
\frac{u_{1}(\alpha)}{v_{2}(\alpha)},\frac{u_{2}(\alpha)}{v_{1}(\alpha)},\frac{u_{2}(\alpha)}{v_{2}(\alpha)}\right\}
$$
\end{enumerate}
\end{definition}

\noindent we adopt the general definition of a fuzzy number given in \cite{rt}

\section{Fuzzy heat-like equations }
In this section, we consider the heat-like equations in one  and tow dimension which can be written in the form
\begin{itemize}
\item One-dimensional :
 \begin{equation}
    U_{t}(t,x) + P(x,\gamma)U_{xx}(t,x) = F(t,x,k) \label{1}
 \end{equation}
\item Two-dimensional :
 \begin{equation}
   U_{t}(t,x,y) + P(x,\gamma)U_{xx}(t,x,y) + Q(y,\beta)U_{yy}(t,x,y) = F(t,x,y,k) \label{x2}
\end{equation} 
or
 \begin{equation}
    U_{t}(t,x,y) + Q(y,\beta)U_{xx}(t,x,y) + P(x,\gamma)U_{yy}(t,x,y) = F(t,x,k) \label{x3}
 \end{equation}
\end{itemize}
subject to certain initial and boundary conditions.

\vskip 0.2in
\noindent These initial and boundary conditions, in state two-dimensional, can come in a variety of forms such as
\begin{center}
$U(0,x,y) = c_1$ or $U(0,x,y) = g_1\left(x,y,c_2\right)$ or $U\left(M_1,x,y\right) = g_2\left(x,y,c_3,c_4\right)$, \ldots 
\end{center}

\noindent In this paper the method is applied for the heat-like equation \eqref{x2}. For \eqref{1} and \eqref{x3}, it is similar to \eqref{x2}, so we will omit them. In following lines, components of \eqref{x2} are enumerated :
\begin{itemize}
\item $I_{j} = [0,M_{j}]$ are three intervals, which  $ M_{j}> 0 \quad  j =1,3 $
\item $F(t,x,y,k)$, $U(t,x,y)$, $P(x,\gamma)$ and $Q(y,\beta)$ will be continuous functions for  $\displaystyle (t,x,y) \in \prod_{j=1}^{3}I_{j}$.
\item $P(x,\gamma)$ and $Q(y,\beta)$ have a finite number of roots for each $(x,y) \in I_{2}\times I_{3}$
\item $k = \left(k_1, \ldots, k_n \right)$, $ c = \left(c_1, \ldots, c_m \right)$, $\gamma = \left(\gamma_1, \ldots, \gamma_s \right)$ and $\beta = \left(\beta_1, \ldots, \beta_e \right)$  are vectors of constants with $k_j\in J_j \subset \R$ , $c_i\in l_i \subset \R$ and $\gamma_r \in H_r\subset \R$, $\beta_l \in D_l \subset \R$. 
\end{itemize}

\noindent Assume that \eqref{x2} has a solution
\begin{equation}
   U(t,x,y) = G(t,x,y,k,c,\gamma,\beta)  \label{2}
\end{equation}
for continuous $G(t,x,y,k,c,\gamma,\beta)$

\noindent and
$$
	\left(G(t,x,y,k,c,\gamma,\beta)\right)_t + P(x,\gamma) \left(G(t,x,y,k,c,\gamma,\beta)\right)_{xx} + Q(y,\beta) \left(G(t,x,y,k,c,\gamma,\beta) \right)_{yy}
$$
is  continuous for  $\displaystyle (t,x,y)\in \prod_{j=1}^{3} I_j$, $\displaystyle k \in J=\prod_{j=1}^{n} J_j$, $\displaystyle c \in L= \prod_{i=1}^{m} L_i$, $\displaystyle \gamma \in H= \prod_{r=1}^{s} H_r$ and $\displaystyle \beta \in D=\prod_{l=1}^{e} D_l$

\noindent Now suppose that $k_j$, $c_i $ ,$\gamma_r$ and $\beta_l$ are imprecise. We will modelling this uncertainty by substitute triangular fuzzy numbers for the  $k_j$, $c_i$ , $\gamma_r$ and $\beta_l$.

\noindent If we fuzzify \eqref{x2}, then we obtain the fuzzy heat-like equation. Using the extension principle we compute  
$\overline{F} \left(t,x,y,\overline{K}\right)$, $\overline{P}\left(x,\overline{\gamma}\right)$ and $\overline{Q}\left(y,\overline{\beta}\right)$ from $F$, $P$ and $Q$ where $\overline{K} = \left(\overline{k}_1, \ldots, \overline{k}_n\right)$, $\overline{\gamma} = \left(\overline{\gamma}_1, \ldots, \overline{\gamma}_s\right)$ and $\overline{\beta} = \left(\overline{\beta}_{1}, \ldots, \overline{\beta}_e\right)$ for $k_j$ , $\gamma_r$ and $\beta_e$
are triangular fuzzy numbers in $J_j$ \ ($0\leq j \leq n$), $H_r$ \ ($0\leq r \leq s$) and $D_l$ \ ($0\leq l\leq e$).

\noindent The function $U$ became $\overline{U}$ where $\overline{U} : \prod_{j=1}^{3}I_{j} \rightarrow \F(\R)$. That is,  $\overline{U}(t,x,y)$ is a fuzzy number. The fuzzy heat-like equation is 
\begin{equation}
   \overline{U}_t(t,x,y) + \overline{P} \left(x,\overline{\gamma}\right) \overline{U}_{xx}(t,x,y) + 
   \overline{Q}\left(y,\overline{\beta}\right) \overline{U}_{yy}(t,x,y)
   = \overline{F} \left(t,x,y,\overline{K}\right)
   \label{3}
\end{equation}
subject to certain initial and boundary conditions.

\noindent The initial and boundary conditions are in the form
\begin{center}
$\overline{U}(0,x,y) = \overline{C}_1$ or $\overline{U}(0,x,y) = \overline{g}_1 \left(x,y,\overline{C}_2\right)$ or  $\overline{U}(M_{1},x,y) = \overline{g}_2 \left(x,y,\overline{C}_3,\overline{C}_4\right)$.
\end{center}
The  $\overline{g}_j $ is the extension principle of $g_j$.

\noindent We wish to solve the problem given in \eqref{3}. Finally, we fuzzify $G$ in \eqref{2}. 

\noindent Let $\overline{Z}(t,x,y) = \overline{G} \left(t,x,y,\overline{K},\overline{C},\overline{\gamma},\overline{\beta}\right)$
where $\overline{Z}$ is computed using the extension principale and is a fuzzy solution. In section 5, we will discuss solution with the same strategy as Buckley-Feuring for fuzzy heat-like equation with
\begin{center}
	$\displaystyle \overline{K}[\alpha] = \prod_{j=1}^{n}\overline{K}_j [\alpha]$,
	$\displaystyle \overline{\gamma}[\alpha] = \prod_{r=1}^{s}\overline{\gamma}_{r}[\alpha]$,
	$\displaystyle \overline{C}[\alpha] = \prod_{i=1}^{m}\overline{C}_{i}[\alpha]$ and 
	$\displaystyle \overline{\beta}[\alpha] = \prod_{l=1}^{e}\overline{\beta}_{l}[\alpha]$
\end{center}

\section{The variational iteration method}
To illustrate the basic idea of the VIM we consider the following model PDE
\begin{equation*}
   L_t\,U + L_x\,U + L_y\,U + N\,U = F(t,x,y,k)
\end{equation*}
where $L_t$, $L_x$ and $L_y$ are linear operators of $t$, $x$ and $y$, respectively, and $N$ is a nonlinear operator, also $F(t,x,y,k)$ is the source non-homogeneous term. According to the VIM, we can express the following correction function in $t$-direction as follows
\begin{equation*}
   U_{n+1}(t,x,y) = U_n(t,x,y) + \int_0^t \lambda \{ L_s\,U_n + \left(L_x + L_y + N\right) \widetilde{U}_n - F(s,x,y,k) \} ds
\end{equation*}
where $\lambda$ is general lagrange multiplier \cite{ttta}, which can be identifier optimally via the variational theory
\cite{tttt,in}, and  $\widetilde{U}_n$ is a restricted variation which means  $\delta \widetilde{U}_n = 0$.

\noindent By this method, we determine first the lagrange multiplier $\lambda$ which will be identified optimally. The approximations 
$U_{n+1}$,  $n \geq 0 $, of the solution $U$ will be readily obtained by suitable choice of trial function $U_0$. Consequently, the solution is given as
\begin{equation*}
 U(t,x,y) = \lim_{n\rightarrow \infty} U_n (t,x,y)
\end{equation*}
According to the  VIM, we construct a correction functional for \eqref{x2} in the form
\begin{equation}
  U_{n+1}(t,x,y) = U_n(t,x,y) + \int_0^t\lambda(s) \{(U_n)_s + P(x,\gamma)(\widetilde{U}_n)_{xx} + Q(y,\beta)(\widetilde{U}_n)_{yy}- F(s,x,y,k) \} ds    \label{4}
\end{equation}
where $n \geq 0$ and $\lambda$ is a lagrange multiplier. Making \eqref{4} stationary with respect to $U_n$, we have 
$$
    \lambda^{'}(s) = 0 \ , \quad  
    1 + \lambda(s)|_{s=t} = 0 
$$
hence, the lagrange multiplier is  $\lambda = -1$. Submitting the results into \eqref{4} leads to the following iteration formula
\begin{equation}
 U_{n+1}(t,x,y) = U_n(t,x,y) - \int_0^t \{(U_n)_s + P(x,\gamma)(\widetilde{U}_n)_{xx} + Q(y,\beta)(\widetilde{U}_n)_{yy} - F(s,x,y,k) \}ds
 \label{o}
\end{equation}
Iteration formula start with initial approximation, for example $U_0(t,x,y) = U(0,x,y)$. Also the VIM used for system of linear and nonlinear partial differential equation  \cite{in} which handled in obtain  Seikkala solution.

\section{Buckley-Feuring Solution (BFS) and  Seikkala solution (SS)}
\subsection{Buckley-Feuring solution}
Buckley-Feuring first present the BFS \cite{ed,tf}.They define for all $t$, $x$, $y$ and $\alpha$,
\begin{equation*}
  \overline{Z}(t,x,y)[\alpha]=[z_{1}(t,x,y,\alpha),z_{2}(t,x,y,\alpha)] \ , \quad 
 \overline{F}(t,x,y,\overline{k})[\alpha] = [F_{1}(t,x,y,\alpha),F_{2}(t,x,y,\alpha)] \ , \quad \forall \alpha \in [0,1]
\end{equation*}
and to check \eqref{3} we must compute $\overline{P}(x,\overline{\gamma})$ and $\overline{Q}(y,\overline{\beta})$. $\alpha$-cuts of $\overline{P}(x,\overline{\gamma})$ and $\overline{Q}(y,\overline{\beta})$ can be found as follows:
\begin{equation*}
 \overline{P}(x,\overline{\gamma})[\alpha]=[P_{1}(x,\alpha),P_{2}(x,\alpha)] \ , \quad \overline{Q}(y,\overline{\beta})[\alpha]=[Q_{1}(y,\alpha),Q_{2}(y,\alpha)] \ , \quad \forall \alpha \in[0,1]
\end{equation*}
that by definition
\begin{equation}
   z_{1}(t,x,y,\alpha)=\min\left\{G(t,x,y,k,c,\gamma,\beta)|k \in \overline{K}[\alpha],c \in \overline{C}[\alpha], \ \gamma \in \overline{\gamma}[\alpha]\  and \ \beta \in \overline{\beta}[\alpha] \right\}
\label{z1}
\end{equation}
\begin{equation}
  z_{2}(t,x,y,\alpha)=\max\left\{G(t,x,k,y,c,\gamma,\beta)|k \in \overline{K}[\alpha],c \in \overline{C}[\alpha],\  \gamma \in \overline{\gamma}[\alpha]\  and \ \beta \in \overline{\beta}[\alpha]\right\}
\label{z2}
\end{equation}
$\displaystyle \forall (t,x,y) \in \prod_{j=1}^{3} I_{j} \quad  \text{and} \quad  \alpha \in [0,1]$

\noindent and
\begin{equation}
  F_{1}(t,x,y,\alpha) = \min\left\{F(t,x,y,k)| k \in \overline{K}[\alpha] \right\}   \label{f1}
\end{equation}
\begin{equation}
    F_{2}(t,x,y,\alpha)=\max\left\{F(t,x,y,k)| k \in \overline{K}[\alpha] \right\}  \label{f2}
\end{equation}
$\displaystyle \forall (t,x,y) \in \prod_{j=1}^{3}I_{j} \  and \ \alpha \in [0,1] $

\noindent and
\begin{equation}
  P_{1}(x,\alpha)=\min\left\{P(x,\gamma)|  \gamma \in \overline{\gamma}[\alpha] \right\} \label{p1}
\end{equation}
\begin{equation}
    P_{2}(x,\alpha)=\max\left\{P(x,\gamma)| \gamma \in \overline{\gamma}[\alpha]\right\} \label{p2}
\end{equation}
$\forall x \in  I_{2} $ and  $ \alpha \in [0,1] $

\noindent and
\begin{equation}
  Q_{1}(y,\alpha)=\min\left\{Q(y,\beta)|  \beta \in \overline{\beta}[\alpha] \right\}  \label{q11}
\end{equation}
\begin{equation}
    Q_{2}(y,\alpha)=\max\left\{Q(y,\beta)| \beta \in \overline{\beta}[\alpha]\right\}  \label{q12}
\end{equation}
$\forall y \in  I_{3}\  et \ \alpha \in [0,1] $

\noindent Assume that  $P(x,\gamma)>0 \ (P_{1}(x,\alpha)> 0)$, $Q(y,\beta)> 0 \ (Q_{1}(y,\alpha)> 0)$ and the  $ z_{i}(t,x,y,\alpha)$ $i=1,2$, have continuous partial so that
$(z_{i})_{t} + P_{i}(z_{i})_{xx} + Q_{i}(z_{i})_{yy}$ is continuous for all $\displaystyle t,x,y \in \prod_{j=1}^{3} I_{j}$ and all $\alpha \in [0,1]$.

\noindent Define
\begin{equation*}
	\Gamma(t,x,y,\alpha) = \Bigl[ (z_{1})_{t} + P_{1}(x,\alpha)(z_{1})_{xx} + Q_{1}(y,\beta)(z_{1})_{yy},(z_{2})_{t}+P_{2}(x,\alpha)(z_{2})_{xx}+Q_{2}(y,\beta)(z_{2})_{yy} \Bigr]
\end{equation*}
for  all $(t,x,y)\in \prod_{j=1}^{3}I_{j}$ and  all $\alpha $.
If, for each fixed $t,x,y \in \prod_{j=1}^{3}I_{j}$, $ \Gamma(t,x,y,\alpha)$ defines the $\alpha-$cut of a fuzzy number, then will be said that  $\overline{Z}(t,x,y)$ is differentiable and is written
\begin{equation*}
 \overline{Z}_{t}[\alpha]+\overline{P}[\alpha]\overline{Z}_{xx}[\alpha]+\overline{Q}[\alpha]\overline{Z}_{yy}[\alpha] = \Gamma(t,x,y,\alpha)
\end{equation*}
for all $(t,x,y)\in \prod_{j=1}^{3}I_{j} $ and  all $\alpha$

\noindent Sufficient conditions for  $\Gamma(t,x,y,\alpha) $ to define  $\alpha-$cut of a fuzzy number was studied in \cite{rt} :
\begin{itemize}
  \item[(i)] $(z_{1})_{t}(t,x,y,\alpha) + P_1(x,\alpha)(z_{1})_{xx}(t,x,y,\alpha) + Q_1(y,\alpha)(z_{1})_{yy}(t,x,y,\alpha)$ is an increasing function of  $\alpha$ for each  $(t,x,y)\in \prod_{j=1}^{3}I_{j} $
  \item[(ii)] $(z_{2})_{t}(t,x,y,\alpha)+P_{2}(x,\alpha)(z_{2})_{xx}(t,x,y,\alpha)+Q_{2}(y,\alpha)(z_{2})_{yy}(t,x,y,\alpha)$ is an decreasing function of $\alpha$ for each $(t,x,y)\in \prod_{j=1}^{3}I_{j}$ and
  \item[(iii)] For $(t,x,y)\in \prod_{j=1}^{3}I_{j}$
\begin{multline}
(z_{1})_{t}(t,x,y,1)+P_{1}(x,1)(z_{1})_{xx}(t,x,y,1)+Q_{1}(y,1)(z_{1})_{yy}(t,x,y,1) \leq (z_{2})_{t}(t,x,y,1) \\
 + P_{2}(x,1)(z_{2})_{xx}(t,x,y,1)+Q_{2}(y,1)(z_{2})_{yy}(t,x,y,1) 
\end{multline}
\end{itemize}

\noindent Now can assume that the $z_{i}(t,x,y,\alpha)$ have continuous partial so $(z_{i})_{t}+P_{i}(x,\alpha)(z_{i})_{xx}+Q_{i}(y,\alpha)(z_{i})_{yy}$ is continuous on  $ \prod_{j=1}^{3}I_{j} \times [0,1] i=1,2$. Hence, if conditions (i)-(iii) above hold, $\overline{Z}(t,x,y)$ is differentiable.\\
For  $\overline{Z}(t,x,y)$ to be a BFS of the fuzzy heat-like equation we need
\begin{itemize}
   \item[(a)] $ \overline{Z}(t,x,y) $ differentiable
   \item[(b)] Eq.(\ref{3})holds for  $\overline{U}(t,x,y)=\overline{Z}(t,x,y)$, and
   \item[(c)] $ \overline{Z}(t,x,y) $  satisfies the initial and boundary conditions. Since no exist specified any particular initial and boundary conditions, then only is checked if \eqref{3} hold.
\end{itemize}
$\overline{Z}(t,x,y)$ is a BFS (without the initial and boundary conditions) if $\overline{Z}(t,x,y)$ is differentiable and 
$$
	(\overline{Z})_{t} + \overline{P}\left(x,\overline{\gamma}\right) (\overline{Z})_{xx} + \overline{Q}\left(y,\overline{\beta}\right) (\overline{Z})_{yy} = \overline{F}\left(t,x,y,\overline{k}\right)
$$ 
or the following equations must hold
\begin{equation} 
  \left(z_1\right)_t + P_1 \left(x,\alpha\right) \left(z_1\right)_{xx} + Q_1 (y,\alpha) (z_1)_{yy} = F_1(t,x,y,\alpha)
\label{e1}
\end{equation}
\begin{equation}
  \left(z_2\right)_t + P_2 (x,\alpha)(z_2)_{xx} + Q_2 (y,\alpha)(z_2)_{yy} = F_2 (t,x,y,\alpha)
\label{e2}
\end{equation}
for all $(t,x,y)\in \prod_{j=1}^{3}I_{j}$ and $\alpha \in [0,1]$.

\noindent Now we will present a sufficient condition for the BFS to exist such as Buckley and Feuring.

\noindent Since there are such a variety of possible initial and boundary conditions, hence we will omit them from the following theorem. One must separately check out the initial and boundary conditions. So, we will omit the constants $c_{i}$, $1\leq i \leq m $, from the problem. Therfore, \eqref{2} becomes
$U(t,x,y) = G(t,x,y,k,\gamma,\beta)$, so $\overline{Z}(t,x,y) = \overline{G}(t,x,y,\overline{K},\overline{\gamma},\overline{\beta})$.

\begin{theoreme}
Suppose $\overline{Z}(t,x,y)$ is differentiable.
 \begin{itemize}
   \item[(a)] 
   \begin{equation}
     \text{if }\  P\left(x,\gamma_i\right)>0  \text{ and } \ \frac{\partial P}{\partial \gamma_i} \frac{\partial G}{\partial \gamma_i} > 0  
     \quad x \in I_2 \  \text{for } \ \  i =1, 2, \ldots, m   \label{5x}
  \end{equation}
and
  \begin{equation}
    \text{if } Q\left(y,\beta_l\right)>0 \text{ and } \frac{\partial Q}{\partial \beta_l} \frac{\partial G}{\partial \beta_l} > 0 \quad 
    y \in I_3 \text{ for } \  l =1, 2, \ldots, e  \label{x5}
  \end{equation}
and
  \begin{equation}
   \text{if }  \frac{\partial G}{\partial k_j} \frac{\partial F}{\partial k_j} > 0 \text{ for } j = 1, 2, \ldots, n  \label{6x}
  \end{equation}
Then BFS=$\overline{Z}(t,x,y)$
\item[(b)] If relations \eqref{5x} does not hold for some i or relation \eqref{x5} does not hold for some $l$, or relation \eqref{6x} does not hold for some $j$, then $\overline{Z}(t,x,y)$ is not a BFS.
\end{itemize}
\label{te}
\end{theoreme}

\noindent \textbf{Proof :}

(a) \ For simplicity assume  $k_j = k$, $\gamma_i = \gamma$, $\beta_l=\beta$ and $\frac{\partial G}{\partial k} < 0$, $\frac{\partial F}{\partial k} < 0$, $\frac{\partial P}{\partial \gamma} > 0$, $\frac{\partial G}{\partial \gamma} > 0$, $\frac{\partial Q}{\partial \beta} < 0$ and $\frac{\partial G}{\partial \beta} < 0$. 

\noindent The proof for $\frac{\partial G}{\partial k} > 0$, $\frac{\partial F}{\partial k} > 0$, 
$\frac{\partial P}{\partial \gamma} < 0$, $\frac{\partial G}{\partial \gamma} < 0$, $\frac{\partial Q}{\partial \beta} > 0$ and 
$\frac{\partial G}{\partial \beta} > 0$ is similar.

\noindent Since $\frac{\partial G}{\partial k} < 0$, $\frac{\partial G}{\partial \gamma} > 0$ and $\frac{\partial G}{\partial \beta} < 0$. Then from \eqref{z1} and \eqref{z2} we have
\begin{equation*}
   z_1 (t,x,y,\alpha) = G\left(t,x,y,k_{2}(\alpha),\gamma_1 (\alpha), \beta_2(\alpha)\right) \ , \quad 
   z_2(t,x,y,\alpha) = G\left(t,x,y,k_1(\alpha), \gamma_2 (\alpha), \beta_1(\alpha)\right)
\end{equation*}
from \eqref{f1}, \eqref{f2} and  $\frac{\partial F}{\partial k} < 0$ we have
\begin{equation*}
  F_1(t,x,y,\alpha) = F\left(t, x, y, k_2(\alpha)\right) \ , \quad 
  F_2(t,x,y,\alpha) = F\left(t, x, y, k_1(\alpha)\right)
\end{equation*}
since, with \eqref{p1}, \eqref{p2} and $\frac{\partial P}{\partial \gamma} > 0$ we have
\begin{equation*}
  P_1(x,\alpha) = P\left(x, \gamma_1(\alpha)\right) \ , \quad 
  P_2(x,\alpha) = P\left(x, \gamma_2(\alpha)\right)
\end{equation*}
from \eqref{q11}), \eqref{q12} and  $\frac{\partial Q}{\partial \beta} < 0$ we have
\begin{equation*}
  Q_1(y,\alpha) = Q\left(y, \beta_2(\alpha)\right) \ , \quad 
  Q_2(y,\alpha) = Q\left(y, \beta_1(\alpha)\right)
\end{equation*}
for all $\alpha \in [0,1]$ where $\overline{K}[\alpha]=[k_{1}(\alpha),k_{2}(\alpha)]\  \overline{\gamma}[\alpha]=[\gamma_{1}(\alpha),\gamma_{2}(\alpha)]$ and $\overline{\beta}[\alpha]=[\beta_{1}(\alpha),\beta_{2}(\alpha)]$.\\
Now $ G(t,x,y,k,\gamma,\beta)$ solves the Eq. (\ref{x2}), which means
\begin{equation*}
  G_t + P(x, \gamma) G_{xx} + Q(y, \beta) G_{yy} = F(t,x,y,k)
\end{equation*}
for all $(t,x,y) \in \prod_{j=1}^{3}I_{j},\  k \in J $ ,$ \gamma \in H$ and $ \beta \in D$

\noindent Assume $\overline{Z}(t,x,y)$ is differentiable and $P(x,\gamma) > 0 $ and $Q(y,\beta)> 0$ so  
\begin{equation}
  \left(z_1(t,x,y,\alpha)\right)_t + P_1(x,\alpha)\left(z_1(t,x,y,\alpha)\right)_{xx} + Q_1(y,\alpha)\left(z_1(t,x,y,\alpha)\right)_{yy} = F_{1}(t,x,y,\alpha) \label{l1}
\end{equation}
\begin{equation}
  \left(z_2(t,x,y,\alpha)\right)_t + P_2(x,\alpha)\left(z_2(t,x,y,\alpha)\right)_{xx} + Q_2(y,\alpha)\left(z_2(t,x,y,\alpha)\right)_{yy} = F_{2}(t,x,y,\alpha) \label{l2}
\end{equation}
for all $(t,x,y) \in \prod_{j=1}^{3}I_{j} $  and $\alpha \in [0,1]$. Hence, \eqref{e1} and \eqref{e2} holds and $\overline{Z}(t,x,y)$ is a BFS.

\medskip
(b) Now consider the situation where \eqref{5x} or \eqref{x5} or \eqref{6x} does not hold. Let us only look at one case where $\frac{\partial Q}{\partial \beta} < 0$ ( assume $\frac{\partial G}{\partial k} > 0$, $\frac{\partial F}{\partial k} > 0$, $ \frac{\partial G}{\partial \gamma} > 0$, $\frac{\partial P}{\partial \gamma} > 0$ and $\frac{\partial G}{\partial \beta} > 0$, $P(x,\gamma) > 0 $ and $Q(y,\beta) > 0 )$. Then we have 
\begin{equation*}
  z_{1}(t,x,y,\alpha) = G(t,x,y,k_{1}(\alpha),\gamma_{1}(\alpha),\beta_{1}(\alpha)) \ , \quad 
  z_{2}(t,x,y,\alpha) = G(t,x,y,k_{2}(\alpha),\gamma_{2}(\alpha),\beta_{2}(\alpha))
\end{equation*}
\begin{equation*}
  F_{1}(t,x,y,\alpha) = F(t,x,y,k_{1}(\alpha)) \ , \quad
  F_{2}(t,x,y,\alpha) = F(t,x,y,k_{2}(\alpha))
\end{equation*}
and
\begin{equation*}
  P_{1}(x,\alpha)=P(x,\gamma_{1}(\alpha)), \qquad P_{2}(x,\alpha)=P(x,\gamma_{2}(\alpha)) 
\end{equation*}
\begin{equation*}
  Q_{1}(y,\alpha)=Q(y,\beta_{2}(\alpha)), \qquad Q_{2}(y,\alpha)=Q(y,\beta_{1}(\alpha))
\end{equation*}
then we have
\begin{equation*}
  (z_{1}(t,x,y,\alpha))_{t}+P_{1}(x,\alpha)(z_{1}(t,x,y,\alpha))_{xx}+Q_{1}(y,\alpha)(z_{1}(t,x,y,\alpha))_{yy}=F_{1}(t,x,y,\alpha)
\end{equation*}
\begin{equation*}
  (z_{2}(t,x,y,\alpha))_{t}+P_{2}(x,\alpha)(z_{2}(t,x,\alpha))_{xx}+Q_{2}(y,\alpha)(z_{2}(t,x,y,\alpha))_{yy}=F_{2}(t,x,y,\alpha)
\end{equation*} 
which is not true. because
\begin{multline}
\left(G\left(t,x,y,k_1(\alpha),\gamma_1(\alpha),\beta_1(\alpha)\right)\right)_t + P\left(x,\gamma_1(\alpha)\right)\left(G\left(t,x,y,k_1(\alpha),\gamma_1(\alpha),\beta_1(\alpha)\right)\right)_{xx} \\
+ Q\left(x,\beta_2(\alpha)\right)\left(G\left(t,x,y,k_1(\alpha),\gamma_1(\alpha),\beta_1(\alpha)\right)\right)_{yy} = F\left(t,x,y,k_1(\alpha)\right)
\end{multline}
\begin{multline}
\left(G\left(t,x,y,k_2(\alpha),\gamma_2(\alpha),\beta_2(\alpha)\right)\right)_{t} + P\left(x,\gamma_1(\alpha)\right)\left(G\left(t,x,y,k_2(\alpha),\gamma_2(\alpha),\beta_2(\alpha)\right)\right)_{xx} \\
+ Q\left(y,\beta_1(\alpha)\right)\left(G\left(t,x,k_1(\alpha),\gamma_1(\alpha),\beta_2(\alpha)\right)\right)_{yy} = F\left(t,x,y,k_2(\alpha)\right)
\end{multline}
\qed

Therefore, if $\overline{Z}(t,x,y) $ is a BFS and it satisfies the initial and boundary conditions we will say that $\overline{Z}(t,x,y) $ is a BFS satisfying the initial and boundary conditions. If $\overline{Z}(t,x,y) $ is not a BFS, then we will consider the SS.

\subsubsection{Seikkala solution (SS)}
Now let us define the SS \cite{rnn}, so
$$ 
	\overline{U}(t,x,y)[\alpha] = \Bigl[ u_{1}(t,x,y,\alpha),u_{2}(t,x,y,\alpha) \Bigr]
$$
For example suppose  $ P(x,\gamma) < 0$ and $ Q(y,\beta) > 0$, so consider the system of heat-like equations
\begin{equation}
  (u_{1})_{t}+P_{1}(x,\alpha)(u_{2})_{xx}+Q_{1}(y,\alpha)(u_{1})_{yy}=F_{1}(t,x,y,\alpha)  \label{7}
\end{equation}
\begin{equation}
  (u_{2})_{t}+P_{2}(x,\alpha)(u_{1})_{xx}+Q_{2}(y,\alpha)(u_{2})_{yy}=F_{2}(t,x,y,\alpha)   \label{8}
\end{equation}
or if $P(x,\gamma) > 0$, $Q(y,\beta)>0$, $\frac{\partial P}{\partial \gamma} > 0$, $\frac{\partial G}{\partial \gamma} < 0$,
$\frac{\partial Q}{\partial \beta} > 0$, $\frac{\partial G}{\partial \beta} > 0$
\begin{equation}
  (u_1)_{t} + P_1(x,\alpha)(u_1)_{xx} + Q_1(y,\alpha)(u_1)_{yy} = F_1(t,x,y,\alpha)  \label{a}
\end{equation}
\begin{equation}
  (u_2)_t + P_2(x,\alpha)(u_2)_{xx} + Q_2(y,\alpha)(u_2)_{yy} = F_2(t,x,y,\alpha))  \label{b}
\end{equation}
for all $(t,x,y)\in \prod_{j=1}^{3}I_{j}$ and $\alpha \in [0,1]$. 

\noindent We add to \eqref{7} and \eqref{8} any initial and boundary conditions. 

\medskip
\noindent For example, if it was $\overline{U}(0,x,y)=\overline{C}$ then we add
\begin{equation}
  u_{1}(0,x,y,\alpha)=c_{1}(\alpha)  \label{9}
\end{equation}
\begin{equation}
  u_{2}(0,x,y,\alpha)=c_{2}(\alpha)  \label{10}
\end{equation}
where $\overline{C}[\alpha]=[c_{1}(\alpha), c_{2}(\alpha)]$.

\noindent Let $u_{i}(t,x,y,\alpha)$ $i=1,2$ solve \eqref{7} and \eqref{8} with initial and boundary conditions. If
\begin{equation*}
    [ u_{1}(t,x,y,\alpha), u_{2}(t,x,y,\alpha) ]
\end{equation*}
define the $\alpha-$cut of a fuzzy number, for all $ (t,x,y)\in \prod_{j=1}^{3}I_{j}$, then $\overline{U}(t,x,y)$ is the SS.

\noindent We will say that derivative condition holds for fuzzy heat-like equation when \eqref{5x}, \eqref{x5} and \eqref{6x} are true.
\begin{theoreme} \ 

  \begin{enumerate}
    \item If BFS=$\overline{Z}(t,x,y)$, then $ SS=\overline{Z}(t,x,y)$
    \item If $ SS=\overline{Z}(t,x,y)$ and the derivative condition holds, then BFS=$\overline{U}(t,x,y)$
  \end{enumerate}
\end{theoreme}

\noindent \textbf{Proof } : 

\begin{enumerate}
  \item Follows from the definition of BFS and SS.
  \item If SS=$\overline{U}(t,x,y)$ then the Seikkala derivative \cite{tf} exists and since the derivative condition holds, therefore, the following equations are holds
\begin{equation} 
  (u_{1})_{t} + P_{1}(x,\alpha)(u_{1})_{xx} + Q_{1}(y,\alpha)(u_{1})_{yy} = F_{1}(t,x,y,\alpha)
\end{equation}
\begin{equation}
  (u_{2})_{t} + P_{2}(x,\alpha)(u_{2})_{xx} + Q_{2}(y,\alpha)(u_{2})_{yy} = F_{2}(t,x,y,\alpha)
\end{equation}
Also suppose one $k_{j}=k$, $\gamma_{i}=\gamma$, $\beta_{l}=\beta$, $\frac{\partial G}{\partial \gamma} < 0$, $\frac{\partial P}{\partial \gamma} < 0$, $\frac{\partial G}{\partial k} < 0$ and $\frac{\partial F}{\partial k} < 0$, $\frac{\partial G}{\partial \beta} > 0$, $ \frac{\partial Q}{\partial \beta} > 0$ (the other cases are similar and are omitted). We see
\begin{eqnarray*}
  z_{1}(t,x,y,\alpha)&=& G(t,x,y,k_{2}(\alpha),\gamma_{2}(\alpha),\beta_{1}(\alpha)) \\
  z_{2}(t,x,y,\alpha)&=& G(t,x,y,k_{1}(\alpha),\gamma_{1}(\alpha),\beta_{2}(\alpha))\\
  F_{1}(t,x,y,\alpha)&=& F(t,x,y,k_{2}(\alpha)) \\
  F_{2}(t,x,y,\alpha)&=& F(t,x,y,k_{1}(\alpha)) \\
  P_{1}(x,\alpha)&=& P(x,\gamma_{2}(\alpha)) \\
  P_{2}(x,\alpha)&=& P(x,\gamma_{1}(\alpha))\\
  Q_{1}(y,\alpha)&=& Q(y,\beta_{1}(\alpha)) \\
  Q_{2}(y,\alpha)&=& Q(y,\beta_{2}(\alpha))
\end{eqnarray*}
Now look at \eqref{e1}, \eqref{e2} also \eqref{z1} and \eqref{z2}, implies that
\begin{equation*}
  u_{1}(t,x,y,\alpha) = G\left(t, x, y, k_{2}(\alpha), \gamma_{2}(\alpha), \beta_{1}(\alpha)\right) = z_{1}(t, x, y, \alpha)
\end{equation*}
\begin{equation*}
  u_{2}(t,x,y,\alpha) = G\left(t, x, y, k_{1}(\alpha), \gamma_{1}(\alpha), \beta_{2}(\alpha)\right) = z_{2}(t, x, y, \alpha)
\end{equation*}
Therefore $BFS=\overline{U}(t,x,y)$
\end{enumerate}
\qed

\begin{lemme}
Consider \eqref{1}, suppose $\overline{Z}(t,x)$ is differentiable.
  \begin{itemize}
    \item[(a)] 
    \begin{equation}
      \text{if } \ P\left(x,\gamma_i\right)>0 \ \text{ and } \ \frac{\partial P}{\partial \gamma_i} \frac{\partial G}{\partial \gamma_i} > 0 \quad x \in I_2 \  \text{for} \  i = 1, 2, \ldots, m  \label{5}
    \end{equation}and
    \begin{equation}
      \text{if } \ \frac{\partial G}{\partial k_j} \frac{\partial F}{\partial k_{j}} > 0 \  \text{for} \  j = 1, 2, \ldots, n  \label{6}
    \end{equation}
Then BFS=$ \overline{Z}(t,x)$
\item[(b)] If relations \eqref{5} does not hold for some $i$  or relation \eqref{6} does not hold for some $j$, then $\overline{Z}(t,x)$ is not a BFS.
\end{itemize}
\label{to}
\end{lemme}

\textbf{Proof } : 
It is similar to theorem (\ref{te})
\qed

\section{Examples}
We consider the following illustrating examples.

\medskip
\noindent \textbf{Example 1 : }\label{x} We first consider the one-dimensional initial value problem
\begin{equation}
  U_{t} + \frac{\gamma}{2}\,x^{2}\,U_{xx} = k   \label{11}
\end{equation}
subject to the initial condition $U(0,x) = cx^{2}$ and $t \in (0,M_{1}], x \in (0,M_{2}]$. Let $k \in [0,J]$, $\gamma \in [0,H]$ and $c \in [L,0]$ are constants. 

\noindent According to the VIM, a correct functional equation for \eqref{11} from \eqref{o} can be constructed as follows
$$
	U_{n+1}(t,x) = U_n(t,x) - \int_{0}^{t} \{ (U_n)_{s}(s,x) + \frac{\gamma}{2} x^2 (\widetilde{U}_{n})_{xx}(s,x)- F(s,x,k) \}ds
$$
Beginning with an initial approximation $U_0(t,x) = U(0,x) = Cx^2$, we can obtain the following successive approximations
\begin{equation*}
	U_{1}(t,x) = Cx^{2}(1-\gamma t) + kt, \quad U_{2}(t,x) = Cx^{2}(1-\gamma t+\gamma^{2} \frac{t^{2}}{2!}) + kt, \quad 
  U_{3}(t,x) = Cx^{2}(1-\gamma t+\gamma^{2} \frac{t^{2}}{2!}-\gamma^{3} \frac{t^{3}}{3!}) + kt
\end{equation*}
and  $\displaystyle U_{n}(t,x)=Cx^{2}(1-\gamma t+\gamma^{2} \frac{t^{2}}{2!}-\gamma^{3} \frac{t^{3}}{3!}+....+(-1)^{n}\gamma^{n} \frac{t^{n}}{n!})+kt $, $n \geq 1$

\noindent The VIM admits the use of $\displaystyle U(t,x)=\lim_{n\rightarrow \infty}U_{n}(t,x)$, which gives the exact solution
 $$ U(t,x)  =cx^{2}\exp(-\gamma t)  +kt$$
Now we fuzzify $F(t,x,k)$, $ P(x,\gamma)$ and  $G(t,x,k,c,\gamma)= cx^{2}\exp(-\gamma t)+kt $. Clearly
 $\overline{F}(t,x,\overline{K})= \overline{K}$, $\overline{P}(x,\overline{\gamma})= \frac{\overline{\gamma}}{2} x^{2}$ so that   $F_{1}(t,x,\alpha)=k_{1}(\alpha)$, $F_{2}(t,x,\alpha)=k_{2}(\alpha)$, $P_{1}(x, \alpha)= \frac{\gamma_{1}}{2} x^{2}$ and $ P_{2}(x, \alpha)= \frac{\gamma_{2}}{2} x^{2}$. Also  $\overline{G}(t,x,\overline{K},\overline{C},\overline{\gamma}) = \overline{C}x^{2}\exp(-\overline{\gamma} t)+\overline{K}t $, therefore,
 $$ z_{i}(t,x,\alpha)= c_{i}(\alpha)x^{2}\exp(-\gamma_{i}(\alpha) t)+k_{i}(\alpha)t $$ 
For $i=1,2$ et $\overline{C} < 0$ ($\overline{C} =(c_{1},c_{2},c_{3})$ with $c_{3} <0 )$,
 $\overline{K}[\alpha] = [k_{1}(\alpha),k_{2}(\alpha)]$, $\overline{C}[\alpha]=[c_{1}(\alpha),c_{2}(\alpha)]$ and $ \overline{\gamma}[\alpha]=[\gamma_{1}(\alpha),\gamma_{2}(\alpha)]$. 
 
\noindent $\overline{Z}(t,x)$ is differentiable because 
$$ 
	\left(z_{i}(t,x,\alpha)\right)_{t} + \frac{\gamma_{i}(\alpha)}{2}x^{2} \left(z_{i}(t,x,\alpha)\right)_{xx} = k_{i}(\alpha) \qquad \text{for } i=1,2
$$
That is, $(\overline{Z})_{t}+\frac{\overline{\gamma}}{2}x^{2}(\overline{Z})_{xx}=\overline{k}$, a fuzzy number.

\noindent Since $P(x,\gamma) > 0$, $\frac{\partial G}{\partial k} > 0$, $\frac{\partial F}{\partial k} > 0$, $\frac{\partial P}{\partial \gamma} > 0$ et $\frac{\partial G}{\partial \gamma} = -cx^{2}t \exp(-\gamma t)> 0$.  Lemma \ref{to} implies the result that $\overline{Z}(t,x)$  is a BFS. We easily see that $z_{i}(0,x,\alpha)=c_{i}(\alpha) x^{2}$ for  $i=1,2$ , so $\overline{Z}(t,x)$ also satisfies the initial condition. The BFS that satisfies the initial condition may be written as
$$ 
	\overline{Z}(t,x) = \overline{C}x^{2}\,\exp(-\overline{\gamma} t) + \overline{K}t 
$$
for all $(t,x)\in (0,M_{1}]\times (0,M_{2}]$

\medskip
\noindent \textbf{Example 2 : }
Assume $c \in ]0,L]$. $\overline{K}$, $ \overline{\gamma}$ and $\overline{C}$ are triangular fuzzy numbers as in Example \ref{x} with $\overline{C} = \left(c_1, c_2, c_3 \right)$ and $c_{1} > 0$. We need to solve the following system
\begin{equation}
  (u_{1})_{t} + \frac{\gamma_{1}(\alpha)}{2} x^{2}(u_{1})_{xx} = k_{1}(\alpha)  \label{q1}
\end{equation}
\begin{equation}
  (u_{2})_{t} + \frac{\gamma_{1}(\alpha)}{2} x^{2}(u_{2})_{xx} = k_{2}(\alpha)  \label{q2}
\end{equation}
\begin{equation}
  u_{1}(0,x,\alpha) = c_{1}(\alpha)x^{2}  \label{q3}
\end{equation}
\begin{equation}
  u_{2}(0,x,\alpha) = c_{2}(\alpha)x^{2}   \label{q4}
\end{equation}
If the intervals $[u_{1}(t,x,\alpha), u_{2}(t,x,\alpha)]$ define $\alpha$-cuts of a fuzzy number $\overline{U}(t,x)$; then 
SS = $\overline{U}(t,x)$. By VIM, the general solution to Eqs.(\ref{q1})-(\ref{q4}) is
\begin{equation}
  u_{1}(t,x,\alpha) = c_{1}(\alpha)x^{2}\exp(-\gamma_{1}(\alpha) t) + k_{1}(\alpha)t  \label{n1}
\end{equation}
\begin{equation}
   u_{2}(t,x,\alpha) = c_{2}(\alpha)x^{2}\exp(-\gamma_{2}(\alpha) t) + k_{2}(\alpha)t \label{n2}
\end{equation}
Now we denote $\left[u_{1}(t,x,\alpha),u_{2}(t,x,\alpha)\right]$ defines $\alpha-$cut of a fuzzy number on area as  $\mathfrak{R}$

\noindent since $u_{i}(t,x, \alpha)$ are continuous and $u_{1}(t,x,1) = u_{2}(t,x,1)$ then we only require to check if 
$\displaystyle \frac{\partial u_{1}}{\partial \alpha} > 0$ 
and $\displaystyle\frac{\partial u_{2}}{\partial \alpha} < 0$ Since $\overline{K}$, $\overline{C}$ and $\overline{\gamma}$  are triangular fuzzy numbers, hence, we pick simple fuzzy parameter so that $ k'_{1}(\alpha),\  c'_{1}(\alpha)$ and $\gamma'_{1}(\alpha)$  are all positive numbers while  $k'_{2}(\alpha)$, $c'_{2}(\alpha)$ and $\gamma'_{2}(\alpha)$ are negative numbers. The "prime" denotes differentiation with respect to  $\alpha$.

\noindent Then there is a $\lambda > 0$ so that  $k'_{1}(\alpha) = c'_{1}(\alpha) = \gamma'_{1}(\alpha) = \lambda$ and 
$k'_{2}(\alpha) = c'_{2}(\alpha) = \gamma'_{2}(\alpha) = - \lambda$. Hence, for the SS exist we need 
\begin{eqnarray}
  \frac{\partial u_{1}}{\partial \alpha} &=& \lambda\left(x^{2}\exp\left(-\gamma_{1}(\alpha) t\right) - c_{1}(\alpha)tx^{2} \exp\left(-\gamma_{1}(\alpha) t\right) + t\right) > 0 \label{r1} \\
 \frac{\partial u_{2}}{\partial \alpha} &=& -\lambda\left(x^{2}\exp\left(-\gamma_{2}(\alpha) t\right) - c_{2}(\alpha)tx^{2} \exp\left(-\gamma_{2}(\alpha) t\right) + t\right) < 0  \label{r}
\end{eqnarray}
Therefore inequalities (\ref{r1}, \ref{r}) holds if $1-c_{2}(\alpha)t > 0$ for all $\alpha \in [0,1]$.
Hence we may choose $\mathfrak{R} $ by the above assumptions in form as
$$
	\mathfrak{R} = \left\{ (t,x)\, :\, 0 < t \leq \frac{1}{c_3} \ \& \  0 < x \leq M_2 \right\}
$$
and the SS exists on $\mathfrak{R}$ in form (\eqref{n1}, \eqref{n2}).

\medskip
\noindent \textbf{Example 3 : }
Consider the two-dimensional heat-like equation with variable coefficients as
\begin{equation}
	\begin{cases}
		U_{t} + \frac{\gamma}{2} x^{2} U_{xx} + \frac{\beta}{2} y^{2} U_{yy} = kxy & \text{} \\
		U(0,x) = c_{1} y^{2} - c_{2} x^{2} & \text{}
	\end{cases} \label{x11}
\end{equation}
which $t \in (0 , M_{1}]$, $x \in (0 , M_{2})$, $y \in (0,M_{3})$, $k \in [0,J] \ , \gamma \in [0,H]\ , \ c_{1} \in [L,0[$, $ c_{2} \in [0,E]$ and $ \beta \in [0,D]$ \\
Similarly we can establish an iteration formula in the form
\begin{multline}
	U_{n+1}(t,x,y) = U_{n}(t,x,y) - \int_{0}^{t} \left( (U_{n})_{s}(s,x,y) + \frac{\gamma}{2}x^{2}(\widetilde{U}_{n})_{xx}(s,x,y) \right. \\
	+ \frac{\beta}{2}y^{2}(\widetilde{U}_{n})_{yy}(s,x,y)- F(s,x,y,k) \left. \right) ds  \label{d}
\end{multline}
We begin with an initial arbitrary approximation: $U_{0}(t,x,y) = U(0,x,y) = c_{1}y^{2} - c_{2}x^{2}$, and using the iteration formula (\ref{d}), we obtain the following successive approximations 
\begin{gather}
	U_{1}(t,x,y) = c_{1} y^{2} (1 - \beta t) - c_{2} x^{2} (1 - \gamma t) + kxyt \nonumber \\
	U_{2}(t,x,y) = c_{1} y^{2} \left(1 - \beta t + \frac{\beta^2 t^2}{2!}\right) - c_{2} x^{2}\left(1 - \gamma t + \frac{\gamma^2 t^2}{2!}\right) + kxyt \nonumber \\
	U_{3}(t,x,y) = c_{1} y^{2} \left(1 - \beta t + \frac{\beta^2 t^2}{2!} - \frac{\beta^3 t^3}{3!}\right) - c_{2} x^{2} \left(1 - \gamma t + \frac{\gamma^2 t^2}{2!} - \frac{\gamma^3 t^3}{3!}\right) + kxyt \nonumber \\
\intertext{and}
	U_{n}(t,x,y) = c_{1} y^{2} \left(1 - \beta t + \frac{\beta^2 t^2}{2!} + \cdots + (-1)^{n} \frac{\beta^n t^n}{n!}\right) - c_{2} x^2 \left(1 - \beta t + \frac{\beta^2 t^2}{2!} + \cdots + (-1)^{n} \frac{\beta^n t^n}{n!}\right) + kxyt \notag
\end{gather}
Then, the exact solution is given by
$$ 
	U(t,x,y) = c_{1}y^{2} \exp(-\beta t) - c_{2} x^{2} \exp( - \gamma t) + kxyt 
$$
fuzzify $F(t,x,k)$, $P(x,\gamma)$, $Q(y,\beta)$ and 
$\displaystyle G(t,x,k,c,\gamma,\beta) = c_{1} y^{2} \exp(-\beta t) - c_{2} x^{2} \exp(-\gamma t) + kxyt$ producing their $\alpha$-cuts
\begin{gather}
	z_{1}(t,x,y,\alpha) = c_{11} y^{2} \exp(-\beta_{1} t) - c_{22} x^{2} \exp(-\gamma_{1} t) + k_{1} xyt \nonumber \\
	z_{2}(t,x,y,\alpha) = c_{12} y^{2} \exp(-\beta_{2} t) - c_{21} x^{2} \exp(-\gamma_{2} t) + k_{2} x y t \nonumber \\
	F_{1}(t,x,y,\alpha) = k_{1}(\alpha)x y, \qquad F_{2}(t,x,y,\alpha) = k_{2}(\alpha)x y \nonumber \\
	P_{1}(x,\alpha) = \frac{\gamma_{1}(\alpha)}{2} x^{2}, \qquad  P_{2}(x,\alpha) = \frac{\gamma_{2}(\alpha)}{2} x^{2} \nonumber \\
	Q_{1}(x,\alpha) = \frac{\beta_{1}(\alpha)}{2} y^{2}, \qquad Q_{2}(x,\alpha) = \frac{\beta_{2}(\alpha)}{2}y^{2} \notag
\end{gather}
where $\overline{C}_{1} < 0 $,  $\overline{C}_{1} =(c_{11},c_{12},c_{13})$ with $c_{13} < 0$ \\

\noindent and 
\begin{gather}
	\overline{K}[\alpha] = \left[ k_{1}(\alpha),k_{2}(\alpha) \right], \quad \overline{C}_{1}[\alpha] = \left[c_{11}(\alpha),c_{12}(\alpha)\right], \quad \overline{C}_{2}[\alpha] = \left[c_{21}(\alpha),c_{22}(\alpha) \right] \nonumber \\
	\overline{\gamma}[\alpha] = \left[\gamma_{1}(\alpha),\gamma_{2}(\alpha)\right] \ \text{ and } \ \overline{\beta}[\alpha] = \left[\beta_{1}(\alpha),\beta_{2}(\alpha)\right] \notag
\end{gather}
We first check to see if $\overline{Z}(t,x,y)$ is differentiable. We compute
$$
	\left[ (z_{1})_{t} + \frac{\gamma_{1}}{2} x^{2} (z_{1})_{xx} + \frac{\beta_{1}}{2} y^{2} (z_{1})_{yy},  (z_{2})_{t} + \frac{\gamma_{2}}{2} x^{2}(z_{2})_{xx} + \frac{\beta_{2}}{2}y^{2}(z_{2})_{yy} \right]
$$
which are $\alpha-$cuts of $\overline{K}\,xy$ i.e $\alpha-$cuts of a fuzzy number. Hence, $\overline{Z}(t,x,y)$ is differentiable.
Since the partial $F$ and $G$ with respect to $k$, the partial P and G with respect to $\gamma$  and the partial $Q$ and $G$ with respect to $\beta$ and $P(x,\gamma)> 0$, $Q(y,\beta)> 0$ are positive then the theorem (\ref{te}) tells us that  $\overline{Z}(t,x,y)$ is a BFS.

\noindent The initial condition 
$$
	z_{1}(0,x,y) = c_{11}(\alpha) y^{2} - c_{22} x^{2}
$$
$$
	z_{2}(0,x,y) = c_{12}(\alpha) y^{2} - c_{21} x^{2} 
$$ 
which are true. Therefore  $\overline{Z}(t,x,y)$ is a BFS which also satisfies the initial condition. This BFS may be written

\noindent For all $(x,y) \in (0,M_{2}) \times (0,M_{3})$, $t \in (0,M_{1}]$
$$
 \overline{Z}(t,x,y) = \overline{C}_{1} y^{2} \exp(-\overline{\beta} t) - \overline{C}_{1} y^{2} \exp(-\overline{\beta} t) \overline{C}_{2} x^{2} \exp(-\overline{\gamma} t)
$$

\medskip
\noindent \textbf{Example 4 : }
We consider the one-dimensional heat-like model
\begin{equation}
	\begin{cases}
		\left(U(t,x)\right)_{t} + \gamma\left(\frac{1}{2} - x\right) \left(U(t,x)\right)_{xx} = -k x^{2} t^{2} & \text{} \\
		U(0,x) = cx^{2} & \text{}
	\end{cases} \label{12}
\end{equation}
which $\displaystyle t \in (0,1] , x \in \left(0,\frac{1}{2}\right)$, and the value of parameters $k$, $c$ and $\gamma $ are in intervals $[0,J]$, $[0,L]$ and  $[0,H]$, respectively. \\
We can obtain the following iteration formula for the Eq.(\ref{12}) 
\begin{equation}
	U_{n+1}(t,x) = U_{n}(t,x) - \int_{0}^{t} \left( (U_{n})_{s}(s,x) + \gamma(\frac{1}{2} - x) (\widetilde{U}_{n})_{xx}(s,x)+kx^{2}s^{2}\right) ds
\label{13}
\end{equation}
We begin with an initial approximation : $U(0,x) = cx^{2}$. By Eq (\ref{13}), after than two iterations the exact solution is given in the closed forms as \\
$$ 
	U(t,x) = G(t,x,k,c,\gamma) = \frac{\gamma}{12}kt^{4} - \frac{\gamma}{6}kxt^{4} - \frac{1}{3} kx^{2}t^{3} + cx^{2}+2\gamma cxt-\gamma ct 
$$
since $\frac{\partial F}{\partial k }= -x^{2}t^{2} < 0$ and $\displaystyle \frac{\partial G}{\partial k } = -\frac{\gamma x t^{4}}{6} + \frac{\gamma t^{4}}{12} - \frac{t^{3}x^{2}}{3} > 0$ for 
$$
	\text{ For } \quad t\in ]0,1] \quad  \text{ and } \quad  x \in \left]0,\frac{1}{4}(-t+\sqrt{4\gamma t+\gamma^{2}t^{2}})\right[
$$
then there is no BFS (lemma (\ref{to})). We proceed to look for a SS. We must solve 
$$
	(u_{1}(t,x,\alpha))_{t} + \gamma_{1}(\alpha)(\frac{1}{2}-x) (u_{1}(t,x,\alpha))_{xx} = -k_{2}x^{2}t^{2}
$$
$$ 
	(u_{2}(t,x,\alpha))_{t} + \gamma_{2}(\alpha)(\frac{1}{2}-x) (u_{2}(t,x,\alpha))_{xx} = -k_{1}x^{2}t^{2}
$$
subject to  
$$ 
	u_{i}(0,x,\alpha) = c_{i}(\alpha)x^{2}, \qquad \text{ for } i=1,2
$$ 
$\displaystyle \widetilde{k}[\alpha] = \left[k_{1}(\alpha),k_{2}(\alpha)\right]$, $\widetilde{c}[\alpha] = \left[c_{1}(\alpha),c_{2}(\alpha)\right]$ and $\overline{\gamma}[\alpha] = \left[\gamma_{1}(\alpha),\gamma_{2}(\alpha)\right]$. By VIM, the solution is
\begin{equation*}
u_{1}(t,x,\alpha) = \frac{\gamma_{1}(\alpha)}{12} k_{2} (\alpha)t^{4} - \frac{\gamma_{1}(\alpha)}{6} k_{2} (\alpha)xt^{4} - \frac{1}{3}k_{2}(\alpha)x^{2}t^{3} + c_{1}(\alpha) x^{2} + 2\gamma_{1}(\alpha) c_{1}(\alpha)xt - c_{1}(\alpha)\gamma_{1}(\alpha)t
\label{41}
\end{equation*}
$$
u_{2}(t,x,\alpha) = \frac{\gamma_{2}(\alpha)}{12} k_{1} (\alpha)t^{4} - \frac{\gamma_{2}(\alpha)}{6} k_{1} (\alpha)xt^{4} - \frac{1}{3}k_{1}(\alpha)x^{2}t^{3} + c_{2}(\alpha)x^{2} + 2\gamma_{2}(\alpha) c_{2}(\alpha)xt - c_{2}(\alpha)\gamma_{2}(\alpha)t 
$$
Now we denote $\displaystyle \left[ u_{1}(t,x,\alpha),u_{2}(t,x,\alpha) \right]$ defines $\alpha-$cut of a fuzzy number on area as  $\mathfrak{R} $ \\
Since $u_{i}(t,x, \alpha)$ are continuous and $u_{1}(t,x,1) = u_{2}(t,x,1)$ then we only require to check if 
$\displaystyle \frac{\partial u_{1}}{\partial \alpha} > 0$ and $\displaystyle \frac{\partial u_{2}}{\partial \alpha} < 0$. 
Since $\overline{K}$ , $\overline{C}$ and $\overline{\gamma}$  are triangular fuzzy numbers, hence, we pick simple fuzzy parameter so that $k'_{1}(\alpha)= c'_{1}(\alpha)= \gamma'_{1}(\alpha)= \lambda$ and  $k'_{2}(\alpha)= c'_{2}(\alpha) =\gamma'_{2}(\alpha) = - \lambda$. Then, for the SS exist we need   

\begin{multline}
\frac{\partial u_1}{\partial \alpha} =  -\frac{t^4}{12} \left( \lambda k_2(\alpha) + (-\lambda)\gamma_{1}(\alpha) \right)  - \frac{x t^4}{6} \left( -\lambda \gamma_{1}(\alpha) + (\lambda)k_{2}(\alpha) \right) + \frac{x^2 t^3}{3}(-\lambda) \\
 + \lambda x^{2} + 2 \left( \lambda c_{1}(\alpha) + (\lambda)\gamma_{1}(\alpha) \right) xt - \left( \lambda c_{1}(\alpha)+(\lambda)\gamma_{1}(\alpha) \right) t > 0 \nonumber
\end{multline}
\begin{multline}
\frac{\partial u_1}{\partial \alpha} = \lambda \left( -\frac{t^4}{12}\left[ k_2(\alpha)-\gamma_{1}(\alpha)\right] -\frac{x t^4}{6} \left[-\gamma_{1}(\alpha) + k_2(\alpha)\right] \right. \\
\left. - \frac{x^{2}t^{3}}{3} + x^{2} + 2 \left[ c_{1}(\alpha) + \gamma_{1}(\alpha)\right] xt - \left[ c_{1}(\alpha) + \gamma_{1}(\alpha)\right] t \right) > 0 \nonumber
\end{multline}
\begin{multline}
\frac{\partial u_2}{\partial \alpha} =  -\frac{t^4}{12} \left[ -\lambda k_1(\alpha) + (\lambda)\gamma_2(\alpha) \right] -\frac{x t^4}{6} \left[\lambda \gamma_2(\alpha) + (-\lambda)k_1(\alpha) \right] + \frac{x^2 t^3}{3}(\lambda) -\lambda x^2 \\
+ 2\left[-\lambda c_{2}(\alpha) + (-\lambda)\gamma_{2}(\alpha)\right] xt - \left[-\lambda c_{2}(\alpha) + (-\lambda)\gamma_{2}(\alpha)\right] t < 0 \nonumber
\end{multline}
\begin{multline}
\frac{\partial u_{2}}{\partial \alpha} = -\lambda ( -\frac{t^{4}}{12}[ k_{1}(\alpha)-\gamma_{2}(\alpha)] -\frac{xt^{4}}{6}[- \gamma_{2}(\alpha)+k_{1}(\alpha)]-\frac{x^{2}t^{3}}{3} + x^{2} \\
+ 2 [ c_{2}(\alpha) + \gamma_{2}(\alpha)] xt - [ c_{2}(\alpha)+\gamma_{2}(\alpha)] t) < 0 \nonumber
\end{multline}
Therefore inequalities hold if
\begin{equation}
	\frac{t^{4}}{12}[ \overline{\gamma}-\overline{K}] + \frac{x t^4}{6} [ \overline{\gamma} - \overline{K}] - \frac{x^2 t^3}{3} + x^2 + 2[ 	\overline{C} + \overline{\gamma} ] xt - [ \overline{C} + \overline{\gamma}] t > 0  \label{14}
\end{equation}
Let $\overline{S} = \overline{C}+\overline{\gamma}$ and  $\overline{R}= \overline{\gamma}-\overline{K}$ where  $\overline{K}[\alpha]=[k_{1}(\alpha),k_{2}(\alpha)]$, $\overline{\gamma}[\alpha]=[\gamma_{1}(\alpha),\gamma_{2}(\alpha)]$ and $\overline{C}[\alpha]=[c_{1}(\alpha),c_{2}(\alpha)]$.

\noindent For $x \in \left]0,\frac{1}{2}\right[$ and $t \in ]0,1]$, the inequality (\ref{14}) holds if we have
\begin{gather}
	0 < t \leq 1 \quad  , \quad  \frac{g(t,\overline{R},\overline{S})+\sqrt{h(t,\overline{R},\overline{S})}}{l(t)} < x < \frac{1}{2} \nonumber \\
\intertext{with}
	g(t,\overline{R},\overline{S})= -12\overline{S}t-\overline{R}t^{4} \nonumber \\
	h(t,\overline{R},\overline{S}) = 144\overline{S}t + 144\overline{S}^{2} t^2 - 48\overline{S} t^4 - 12\overline{R} t^4 + 24\overline{R} \overline{S} t^{5} + 4\overline{R} t^7 + \overline{R}^{2} t^8 \nonumber
\intertext{and}
	l(t) = 12 - 4\,t^3 \nonumber \\
\intertext{where}
	g\left(t,\overline{R},\overline{S}\right)[\alpha] = \left[g_{1}(t,\alpha), g_{2}(t,\alpha)\right], \quad h\left(t,\overline{R},\overline{S}\right)[\alpha] = \left[h_{1}(t,\alpha),h_{2}(t,\alpha)\right] \notag
\end{gather}
and 
\begin{eqnarray*}
	g_{1}(t,\alpha) &=& \min\{g(t,S,R)|S \in \overline{S}[\alpha], R \in \overline{R}[\alpha]\} \\
  g_{2}(t,\alpha) &=& \max\{g(t,S,R)|S \in \overline{S}[\alpha], R \in \overline{R}[\alpha]\} \\
  h_{1}(t,\alpha) &=& \min\{h(t,S,R)|S \in \overline{S}[\alpha], R \in \overline{R}[\alpha]\}\\
  h_{2}(t,\alpha) &=& \max\{h(t,S,R)|S \in \overline{S}[\alpha], R \in \overline{R}[\alpha]\}
\end{eqnarray*}
\noindent We find that
$$
	\max\left\{ \frac{g(t,\overline{R},\overline{S})[\alpha]+\sqrt{h(t,\overline{R},\overline{S})[\alpha]}}{l(t)} : 0 < t \leq 1 \ \text{ and} \ 0 \leq \alpha \leq 1\right\} = e 
$$
with $e$ is number  in $\R$. Hence we may choose $\mathfrak{R}$ by the above assumptions in form as
$$ 
	\mathfrak{R} = \left\{ (t,x) : 0 < t \leq 1, \quad e \leq x < \frac{1}{2} \right\} 
$$
and the SS exists on $\mathfrak{R}$ in form Eqs.(\ref{41}).

\vskip 0.3in
\noindent \textbf{Example 5 : }
We consider the one-dimensional heat-like model
$$
	\left\{ \begin{array}{ll}
		U_{t}(t,x) - \gamma\, U_{xx} = -k\, \cos(x) & \hbox{} \\
		U(0,x) = c\, \sin(x) & \hbox{}
   \end{array} \right.
$$
which  $\displaystyle x \in (0 , \frac{\Pi}{2})$, $t \in [0,M]$ and the value of parameters $k\in [J,0[$, $c\in [0,L]$ and $\gamma\in (0,R]$. \\
We can obtain the following iteration formula
\begin{equation}
  U_{n+1}(t,x) = U_{n}(t,x) - \int_{0}^{t} \left( (U_{n})_{s}(s,x) - \gamma(\widetilde{U}_{n})_{xx}(s,x) + k \cos(x) \right) ds
\label{100}
\end{equation}
We begin with an initial approximation : $U_{0}(t,x)=U(0,x)=c \sin(x)$. By (\ref{100}), the following successive approximation are obtained 
\begin{gather}
	U_{0}(t,x) = U(0,x) = c \sin(x) \nonumber \\
	U_{1}(t,x) = c \sin(x)(1-\gamma t) - k \cos(x)t \nonumber \\
	U_{2}(t,x) = c \sin(x)\left(1-\gamma t + \frac{\gamma^{2}t^{2}}{2!}\right) + k \cos(x)\left(-t +\frac{\gamma t^{2}}{2!} \right) \nonumber \\
\intertext{and for $n\geq 1$}
	U_{n}(t,x) = c \sin(x) \left( 1-\gamma t + \frac{\gamma^{2}t^{2}}{2!} + \cdots + (-1)^{n}\frac{\gamma^{n}t^{n}}{n!}\right) + \frac{k}{\gamma} \cos(x)\left(-\gamma t + \frac{{\gamma}^2 t^{2}}{2!} + \cdots + (-1)^{n}\frac{\gamma^{n} t^{n}}{n!}\right) \notag
\end{gather}
The VIM admits the use of $ U(t,x) = \lim\limits_{n\rightarrow \infty} U_{n}(t,x)$, which gives the exact solution. There is no BFS because
$$ 
	U(t,x) = G(t,x,k,c,\gamma) = c \sin(x)\exp(-\gamma t ) + \frac{k}{\gamma}\cos(x)\left(\exp(-\gamma t)-1\right)
$$
$P(x,\gamma) = - \gamma <0$ with $\gamma \in (0,R]$ (lemma (\ref{to})). We proceed to look for a SS. 

\noindent We must solve
$$
  \left(u_{1}(t,x,\alpha)\right)_{t}- \gamma_{2} \left(u_{2}(t,x,\alpha)\right)_{xx} = -k_{2} \cos(x)
$$
$$
  \left(u_{1}(t,x,\alpha)\right)_{t} - \gamma_{1} \left(u_{1}(t,x,\alpha)\right)_{xx} = -k_{1} \cos(x)
$$
subject to $u_{i}(0,x,y) = c_{i}(\alpha) \sin(x)$ for $i=1,2$

\noindent and $\displaystyle \overline{K}[\alpha] = \left[k_{1}(\alpha), k_{2}(\alpha)\right] \quad \overline{K} < 0$, 
$\displaystyle \overline{C}[\alpha] = \left[c_{1}(\alpha), c_{2}(\alpha)\right]$ and $\displaystyle \overline{\gamma}[\alpha] = [\gamma_{1}(\alpha), \gamma_{2}(\alpha)]$.

\noindent The solution is
\begin{multline}
	u_{1}(t,x,\alpha) = c_{1}(\alpha)\sin(x)\cosh(\sqrt{\gamma_{1}(\alpha)\gamma_{2}(\alpha)} t) - \frac{c_{2}(\alpha)\gamma_{2}(\alpha)} {\sqrt{\gamma_{1}(\alpha)\gamma_{2}(\alpha)}} \sin(x) \sinh(\sqrt{\gamma_{1}(\alpha)\gamma_{2}(\alpha)} t) \\
	+ \frac{k_{1}(\alpha)}{\gamma_{1}(\alpha)}\cos(x)(\cosh(\sqrt{\gamma_{1}(\alpha)\gamma_{2}(\alpha)} t)-1)-\frac{k_{2}(\alpha)}{\sqrt{\gamma_{1}(\alpha)\gamma_{2}(\alpha)}}\cos(x)\sinh(\sqrt{\gamma_{1}(\alpha)\gamma_{2}(\alpha)} t ) \nonumber
\end{multline}
\begin{multline} 
	u_{2}(t,x,\alpha) = c_{2}(\alpha)\sin(x)\cosh(\sqrt{\gamma_{1}(\alpha)\gamma_{2}(\alpha)} t) - \frac{c_{1}(\alpha)\gamma_{1}(\alpha)} {\sqrt{\gamma_{1}(\alpha)\gamma_{2}(\alpha)}} \sin(x) \sinh(\sqrt{\gamma_{1}(\alpha)\gamma_{2}(\alpha)} t) \\
	+ \frac{k_{2}(\alpha)}{\gamma_{2}(\alpha)}\cos(x)(\cosh(\sqrt{\gamma_{1}(\alpha)\gamma_{2}(\alpha)} t)-1)-\frac{k_{1}(\alpha)}{\sqrt{\gamma_{1}(\alpha)\gamma_{2}(\alpha)}}\cos(x)\sinh(\sqrt{\gamma_{1}(\alpha)\gamma_{2}(\alpha)} t )
\nonumber
\end{multline}
We need only to check if  $\displaystyle \frac{\partial u_{1}}{\partial \alpha} > 0$ and $ \displaystyle \frac{\partial u_{2}}{\partial \alpha} < 0$, since the  $u_{i}(t,x, \alpha)$ are continuous and  $u_{1}(t,x,1) = u_{2}(t,x,1)$.

\noindent  We pick simple fuzzy parameter $k'_1(\alpha) =c'_1(\alpha) = \gamma'_1(\alpha)=\lambda > 0$ and $k'_2(\alpha)= c'_2(\alpha) = \gamma'_2(\alpha)= - \lambda$.

\noindent Let $w=\gamma_{1}(\alpha)\gamma_{2}(\alpha)$. Now we need to check if  $\displaystyle \frac{\partial u_1}{\partial \alpha} > 0$  and $\displaystyle \frac{\partial u_{2}}{\partial \alpha} < 0$, for all $t \in [0,M]$.\\
But this is not true since for large t (we will assume M is a very large positive number )\\
\begin{multline}
	u_{1}(t,x,\alpha)\approx c_{1}(\alpha)\sin(x)\frac{\exp(\sqrt{w} t)}{2} - \frac{c_{2}(\alpha)\gamma_{2}(\alpha)}{\sqrt{w}} \sin(x) \frac{\exp(\sqrt{w} t) }{2} \\
	+ \frac{k_1(\alpha)}{\gamma_1(\alpha)}\cos(x)\frac{\exp(\sqrt{w} t)}{2} - \frac{k_2(\alpha)}{\sqrt{w}} \cos(x) \frac{\exp(\sqrt{w} t )}{2} 
\nonumber
\end{multline}
\begin{multline}
	u_{2}(t,x,\alpha)\approx c_{2}(\alpha)\sin(x)\frac{\exp(\sqrt{w} t)}{2} - \frac{c_{1}(\alpha)\gamma_{1}(\alpha)}{\sqrt{w}} \sin(x) \frac{\exp(\sqrt{w} t ) }{2} \\
	+ \frac{k_{2}(\alpha)}{\gamma_{2}(\alpha)}\cos(x)\frac{\exp(\sqrt{w} t)}{2} - \frac{k_{1}(\alpha)}{\sqrt{w}}\cos(x)\frac{\exp(\sqrt{w} t)}{2} 
\nonumber
\end{multline}
\begin{eqnarray*}
  u_{1}(t,x,\alpha) &\approx& ( A_{1}(x,\alpha))\frac{\exp(\sqrt{w} t )}{2} \\
  u_{2}(t,x,\alpha) &\approx& (A_{2}(x,\alpha))\frac{\exp(\sqrt{w} t )}{2}
\end{eqnarray*}
where
\begin{eqnarray*}
A_{1}(x,\alpha) &=& (c_{1}(\alpha) -\frac{c_{2}(\alpha)\gamma_{2}(\alpha)}{\sqrt{w}})\sin(x)+(\frac{k_{1}(\alpha)}{\gamma_{1}(\alpha)} -\frac{k_{2}(\alpha)}{\sqrt{w}})\cos(x), \\
A_{2}(x,\alpha) &=& (c_{2}(\alpha) -\frac{c_{1}(\alpha)\gamma_{1}(\alpha)}{\sqrt{w}})\sin(x)
+(\frac{k_{2}(\alpha)}{\gamma_{2}(\alpha)} -\frac{k_{1}(\alpha)}{\sqrt{w}})\cos(x)
\end{eqnarray*}
and $A_{1}(x,\alpha) \leq A_{2}(x,\alpha)$ for all $\alpha$, $x \in (0,\frac{\pi}{2}).$
\begin{multline}
 \frac{\partial u_{1}}{\partial \alpha} \approx \lambda( \sin(x)(1+ \frac{2(\gamma_{2}(\alpha) + c_{2}(\alpha))w + (\gamma_{2}(\alpha)-\gamma_{1}(\alpha))\gamma_{2}(\alpha)c_{2}(\alpha)}{2w\sqrt{w}}) \\
 + \cos(x)(\frac{\gamma_{1}(\alpha)-k_{1}(\alpha)}{\gamma_{1}(\alpha)^{2}}  + \frac{(2w+\gamma_{2}(\alpha) - \gamma_{1}(\alpha))k_{2}(\alpha)}{2w\sqrt{w}})+A_{1}(x,\alpha) t \frac{\gamma_{2}(\alpha) - \gamma_{1}(\alpha)}{2\sqrt{w}})\frac{\exp(\sqrt{w}t)}{2}
\nonumber
\end{multline}
\begin{multline}
 \frac{\partial u_{2}}{\partial \alpha} \approx -\lambda( \sin(x)(1+\frac{2(\gamma_{1}(\alpha)+c_{1}(\alpha))w + (\gamma_{2}(\alpha) - \gamma_{1}(\alpha))\gamma_{1}(\alpha)c_{1}(\alpha)}{2w\sqrt{w}}) \\
 + \cos(x)(\frac{\gamma_{2}(\alpha)-k_{2}(\alpha)}{\gamma_{2}(\alpha)^{2}} + \frac{2w +(\gamma_{2}(\alpha) - \gamma_{1} (\alpha))k_{1}(\alpha)}{2w\sqrt{w}}) - A_{2}(x,\alpha) t \frac{\gamma_{2}(\alpha)-\gamma_{1}(\alpha)}{2\sqrt{w}})\frac{\exp(\sqrt{w}t)}{2}
\nonumber
\end{multline}
so as $t$ grows we get $\displaystyle \frac{\partial u_{1}}{\partial \alpha} > 0$ if $A_{1}(x,\alpha)>0$ and $\displaystyle  \frac{\partial u_{2}}{\partial \alpha} < 0$ if $A_{2}(x,\alpha)<0$ this is not true because $A_{1}(x,\alpha) \leq A_{2}(x,\alpha)$ for all $\alpha$ and $x \in (0,\frac{\pi}{2})$

\noindent So, in general SS does not exist for large values of $t$. SS exists for $0 \leq t \leq M$ for some $M > 0$.

\noindent therefore, $\overline{U}(t,x)$ is SS and
\begin{multline}
 \overline{U}(t,x) = \overline{C} \sin(x)\cosh(\sqrt{w}t)-\frac{\overline{C} \overline{\gamma}}{\sqrt{w}} \sin(x) \sinh(\sqrt{w} t) \\
	+ \frac{\overline{K}}{\overline{\gamma}} \cos(x) (\cosh(\sqrt{w}t)-1) -\frac{\overline{K}}{\sqrt{w}}\cos(x)\sinh(\sqrt{w} t)
\nonumber
\end{multline}
for all  $t \in [0,M]$, $x \in (0,\frac{\Pi}{2})$

\section{Conclusion}
In this paper, the same strategy as Buckley-Feuring using the VIM has been successfully applied for solving fuzzy heat-like equations in one and two dimensions with variable coefficients.

\noindent Working procedure is to add to this strategy others functions with fuzzy parameters and by help of the VIM ,we give a prolongement for the strategy of Buckley-Feuring for the proposed models. \\
Application of VIM is easy and calculation of successive approximations is direct and straightforward.\\
In exact fuzzy solution, if the BFS fails to exist we check if the SS exists and when the SS fails to exist we offer no solution to the
fuzzy heat-like equations.


\end{document}